\documentclass[12pt]{amsart} 
\usepackage{amssymb}
\usepackage{fullpage}
\usepackage[colorlinks=true]{hyperref}
\usepackage{graphicx}

\newcommand{\bbC}{\mathbb{C}}
\newcommand{\bbK}{\mathbb{K}}

\newcommand{\bbQ}{\mathbb{Q}}

\newcommand{\bbZ}{\mathbb{Z}}

\newcommand{\cC}{\mathcal{C}}
\newcommand{\cD}{\mathcal{D}}
\newcommand{\cN}{\mathcal{N}}
\newcommand{\cO}{\mathcal{O}}

\DeclareMathOperator{\core}{sf}

\DeclareMathOperator{\rad}{rad}

\DeclareMathOperator{\Tr}{Tr}

\numberwithin{equation}{section}

\newtheorem{theorem}{Theorem}[section]
\newtheorem{conjecture}[theorem]{Conjecture}
\newtheorem{corollary}[theorem]{Corollary}
\newtheorem{lemma}[theorem]{Lemma}
\newtheorem{proposition}[theorem]{Proposition}

\theoremstyle{definition}

\newtheorem*{remark-nonum}{Remark}
\newtheorem*{remarks-nonum}{Remarks}

\newcommand{\lemDio}{2.1}
\newcommand{\lemHypg}{2.2}

\newcommand{\eqGDefn}{2.4}
\newcommand{\eqQDefn}{2.7}
\newcommand{\eqKUB}{2.8}
\newcommand{\eqEDefn}{2.9}

\newcommand{\lemLB}{3.5}
\newcommand{\lemUB}{3.7}
\newcommand{\lemGap}{3.8}
\newcommand{\lemTheta}{3.9}
\newcommand{\lemGcd}{3.10}

\newcommand{\lemGNd}{3.12}

\newcommand{\propRep}{3.1}

\newcommand{\subSectConstr}{2.1} 
\newcommand{\subSectStepiv}{4.5} 

\newcommand{\eqQAbs}{3.25} 
\newcommand{\eqOmegaMinusStuffEquality}{3.30} 
\newcommand{\eqOmegaUB}{3.31} 

\newcommand{\eqELB}{4.5} 
\newcommand{\eqQUB}{4.7} 
\newcommand{\eqYLUBUsingYK}{4.15} 

\begin{document}

\title[Bounds on the number of squares \ldots]{Bounds on the number of squares in recurrence sequences: $y_{0}=b^{2}$ (I)}

\author{Paul M Voutier}
\address{London, UK}
\email{Paul.Voutier@gmail.com}

\date{}

\begin{abstract}
We continue and generalise our earlier investigations of the number of squares in
binary recurrence sequences. Here we consider sequences, $\left( y_{k} \right)_{k=-\infty}^{\infty}$,
arising from the solutions
of generalised negative Pell equations, $X^{2}-dY^{2}=c$, where $-c$ and $y_{0}$
are any positive squares.
We show that there are at most $2$ distinct squares larger than an explicit
lower bound in such sequences. From this result, we also show that there are
at most $5$ distinct squares when $y_{0}=b^{2}$ for infinitely many values of
$b$, including all $1 \leq b \leq 24$, as well as once $d$ exceeds an explicit
lower bound, without any conditions on the size of such squares.
\end{abstract}


\maketitle

\section{Introduction}
\label{sect:intro}

The study of the arithmetic properties of recurrence sequences is an important
part of number theory and has a long history (see \cite{Everest}). Questions
regarding the squares in
such sequences, in particular, binary recurrence sequences, are important in
their own right, as well as for their connection to Diophantine equations of the
form $aX^{2} - bY^{4} = c$. Such equations are also quartic models of elliptic
curves, adding to their significance.

In recent work \cite{V5,V6}, we developed a technique for bounding the number
of distinct squares in binary recurrence sequences. In these papers, we applied
our technique to sequences $\left( y_{k} \right)_{k=-\infty}^{\infty}$ with
$y_{0}=1$ that arise from the solutions of generalised Pell equations,
$X^{2}-dY^{2}=c$, with $-c$ a positive square \cite{V5} or $1$, $2$ or $4$ times
a prime power \cite{V6}. We were able to obtain best possible results
for most of the sequences we considered. Here we extend such results to
when $y_{0}$ is any positive square and for $-c$ a positive square.

\subsection{Notation}
\label{subsect:notation}

Let $a$, $b$ and $d$ be positive integers such
that $d$ is not a square. Suppose $\alpha=a+b^{2}\sqrt{d}$ has $N_{\alpha}
=N_{\bbQ \left( \sqrt{d} \right)/\bbQ}(\alpha)=a^{2}-b^{4}d$ and let
$\varepsilon=\left( t+u\sqrt{d} \right)/2$ be a unit in
$\cO_{\bbQ \left( \sqrt{d} \right)}$ with $t$ and $u$ positive integers.

We define the two sequences $\left( x_{k} \right)_{k=-\infty}^{\infty}$
and $\left( y_{k} \right)_{k=-\infty}^{\infty}$ by
\begin{equation}
\label{eq:yk-defn}
x_{k}+y_{k}\sqrt{d}
=\alpha \varepsilon^{2k}.
\end{equation}

Observe that $x_{0}=a$,
$y_{0}=b^{2}$,
\begin{equation}
\label{eq:yPM1}
y_{1}= \frac{b^{2}\left( t^{2}+du^{2} \right)+2atu}{4}, \hspace*{3,0mm}
y_{-1}=\frac{b^{2}\left( t^{2}+du^{2} \right)-2atu}{4}
\end{equation}
and that both sequences satisfy
the recurrence relation
\begin{equation}
\label{eq:yk-recurrence}
u_{k+1}= \frac{t^{2}+du^{2}}{2} u_{k}-u_{k-1},
\end{equation}
for all $k \in \bbZ$. Note that $\left( t^{2}+du^{2} \right)/2
=\Tr_{\bbQ \left( \sqrt{d} \right)/\bbQ} \left( \varepsilon^{2} \right)$.

Observe that from \eqref{eq:yk-defn},
\[
x_{k}^{2}-dy_{k}^{2}=N_{\alpha}.
\]
So this pair of sequences provides a family of solutions of the generalised
Pell equation, $x^{2}-dy^{2}=N_{\alpha}$.

We are interested here in squares in the sequence of $y_{k}$'s, so we set the
coefficient of $\sqrt{d}$ in $\alpha$ to be a square. Furthermore, we choose
$\alpha$ such that $b^{2}$ is the smallest square among the $y_{k}$'s and let
$K$ be the largest negative integer such that $y_{K}>b^{2}$.

For any non-zero integer, $n$, let $\core(n)$ be the unique squarefree integer
such that $n/\core(n)$ is a square. We will put $\core(1)=1$.

\subsection{Conjectures}
\label{subsect:conj}

We start with some conjectures regarding squares in the sequence of $y_{k}$'s
first stated in \cite{V5}.
The dependence on the arithmetic of $N_{\alpha}$ is noteworthy.

\begin{conjecture}
\label{conj:1-seq}
There are at most four distinct integer squares among the $y_{k}$'s.

If $\core \left( \left| N_{\alpha} \right| \right)|(2p)$ where $p$ is an odd
prime, then there are at most three distinct integer squares among the $y_{k}$'s.

Furthermore, if $\left| N_{\alpha} \right|$ is a perfect square, then there are
at most two distinct integer squares among the $y_{k}$'s.
\end{conjecture}

In fact, a more general result than Conjecture~\ref{conj:1-seq} also appears to
be true. Removing the restriction to even powers of $\varepsilon$ in \eqref{eq:yk-defn},
define $\left( x_{k}' \right)_{k=-\infty}^{\infty}$ and $\left( y_{k}' \right)_{k=-\infty}^{\infty}$
by
\begin{equation}
\label{eq:ykPrime-defn}
x_{k}'+y_{k}'\sqrt{d}
=\alpha \varepsilon^{k}.
\end{equation}

\begin{conjecture}
\label{conj:2-seq}
There are at most four distinct integer squares among the $y_{k}'$'s.

If $\left| N_{\alpha} \right|$ is a prime power or a perfect square,
then there are at most three distinct integer squares among the $y_{k}'$'s.
\end{conjecture}

Computational evidence for these conjectures was presented in Subsection~1.3 of \cite{V5}.

\begin{remark-nonum}
The distinctness condition in these conjectures, and in our results below, is
important, as such sequences can have repeated elements. E.g., $(a,b,d,t,u)=(42,4,7,16,6)$
where $y_{-k}=y_{k-1}$ for all $k \geq 1$.
But this can only happen when $\alpha$ divided by its algebraic conjugate is a
unit in the ring of integers.
\end{remark-nonum}

\subsection{Results}
\label{subsect:results}

In this paper, we obtain a small upper bound on the number of distinct squares
in such sequences when $-N_{\alpha}$ is a square.

\begin{theorem}
\label{thm:1.2-seq-new}
Let $a$, $b$ and $d$ be positive integers, where $d$ is not a square, $N_{\alpha}<0$
and $-N_{\alpha}$ is a square. There are at most two distinct squares among the
$y_{k}$'s with $k \geq 1$ or $k \leq K$, and
\[
y_{k}> \frac{16.33 b^{8/3}\left| N_{\alpha} \right|^{2}}{\sqrt{d}}.
\]
\end{theorem}

\begin{remark-nonum}
This is a somewhat weaker generalisation of Theorem~1.4 in \cite{V5}, where we
proved Conjecture~\ref{conj:1-seq} when $b=1$ and $-N_{a}$ is a positive square
for all but a very specific (albeit infinite) set of sequences. For these remaining
sequences, we showed there are at most $3$ distinct squares.

We state a more precise version of Theorem~\ref{thm:1.2-seq-new} in Proposition~\ref{prop:1.2-seq-new}
in Subsection~\ref{subsect:thm12-consolidation} below.
\end{remark-nonum}

\begin{remark-nonum}
At the expense of supposing there are more distinct squares, we can use our
gap principle in Lemma~\ref{lem:gap} below to reduce the bound on $y_{k}$. For
example, if we suppose there are at most three distinct squares among the
$y_{k}$'s with $k \geq 1$ or $k \leq K$, then we need only also assume that
$y_{k}> 0.66b^{2+2/9}\left| N_{\alpha} \right|^{1+1/3}/d^{1-1/6}$. This is
because, applying Lemma~\ref{lem:gap}, we have $57.32d^{2}/ \left( b^{4}\left| N_{\alpha} \right|^{2} \right)$
times the cube of this lower bound is
$16.47\ldots b^{8/3}\left| N_{\alpha} \right|^{2}/\sqrt{d}$, exceeding the lower
bound in Theorem~\ref{thm:1.2-seq-new}.

As we increase the number of additional distinct squares added, the lower bound
for $y_{k}$ approaches $O \left( b^{2}\left| N_{\alpha} \right|/d \right)$,
where the $O$ constant approaches $0$.
\end{remark-nonum}

\begin{theorem}
\label{thm:1.3-seq-new}
Let $a$, $b$ and $d$ be positive integers, where $d$ is not a square, $N_{\alpha}<0$
and $-N_{\alpha}$ is a square.

{\rm (a)} If $b$ is of the form $b=b_{1}b_{2}$ where $b_{1} \in \left\{ 1,5,13, 17 \right\}$
and $b_{2}$ has no prime factors congruent to $1 \bmod{4}$, then there are at
most two distinct squares among the $y_{k}$'s with $k \geq 2$ or $k \leq K-1$.

{\rm (b)} If $b$ is not of the form in part~{\rm (a)}, then there are at most two
distinct squares among the $y_{k}$'s with $k \geq 2$ or $k \leq K-1$, provided that
\[
d \geq \frac{17 \left| N_{\alpha} \right|^{1/2}b^{28/13}}{u^{24/13}}.
\]
\end{theorem}

The following corollary is immediate from Theorem~\ref{thm:1.3-seq-new}.

\begin{corollary}
\label{cor:1.3-seq-new}
Let $a$, $b$ and $d$ be positive integers, where $d$ is not a square, $N_{\alpha}<0$
and $-N_{\alpha}$ is a square.

{\rm (a)}  If $b$ is of the form $b=b_{1}b_{2}$ where $b_{1} \in \left\{ 1,5,13, 17 \right\}$
and $b_{2}$ has no prime factors congruent to $1 \bmod{4}$, there are at most
five distinct squares among the $y_{k}$'s.

{\rm (b)} If $b$ is not of the form in part~{\rm (a)}, there are at most five
distinct squares among the $y_{k}$'s, provided that
\[
d \geq \frac{17 \left| N_{\alpha} \right|^{1/2}b^{28/13}}{u^{24/13}}.
\]
\end{corollary}

From part~(a) of this corollary, we obtain the following result regarding distinct
squares in the more general sequence $\left( y_{k}' \right)_{k=-\infty}^{\infty}$
defined in \eqref{eq:ykPrime-defn}.

\begin{corollary}
\label{cor:1.4-seq-new}
Let $a$, $b$ and $d$ be positive integers, where $d$ is not a square, $N_{\alpha}<0$
and $-N_{\alpha}$ is a square.

If the smallest squares among the even-indexed $y_{k}'$ and among the odd-indexed
$y_{k}'$ can each be written as $\left( b_{1}b_{2} \right)^{2}$
where $b_{1} \in \left\{ 1,5,13, 17 \right\}$
and $b_{2}$ has no prime factors congruent to $1 \bmod{4}$, then there
are at most $10$ distinct squares among the $y_{k}'$'s.
\end{corollary}

\begin{remark-nonum}
The condition on $b$ in Theorem~\ref{thm:1.3-seq-new}(a), Corollary~\ref{cor:1.3-seq-new}(a)
and Corollary~\ref{cor:1.4-seq-new} is satisfied for all $b<25$.
\end{remark-nonum}

\subsection{Our method of proof}

We use the same approach here as we developed and used in \cite{V5,V6}.
The main differences in the proof here and in \cite{V5,V6} are the following:\\
(1) in \cite{V5,V6}, we only needed $y_{k} \geq 4$. Here we use the lower bound
for $y_{k}$ in Lemma~\lemGap{} of \cite{V5}. In fact, our use of Lemma~\lemGap{}
of \cite{V5} here and in subsequent papers is one of the reasons for the form
of the lower bound for $y_{k}$ there. It means that the $y_{k}$'s are sufficiently
large to overcome the effects of the denominator in the gap principle there;\\
(2) in Section~\ref{sect:proof-thm13}, use $y_{k}$'s with $k$ further from $0$
in our proof to obtain a better lower bounds for the $y_{k}$'s from Lemma~\ref{lem:Y-LB1};\\
(3) assume the existence of an additional square so we can use the gap principle
twice in Subsections~\ref{subsect:thm12-step-i} and \ref{subsect:thm12-step-ii}
to get larger lower bounds for $y_{k_{3}}$;\\
(4) obtaining bounds for $d$ enabling us to treat small $b$. See
Section~\ref{sect:proof-thm13} and Subsection~\ref{subsect:proof-thm13a}.

Lastly, the code used in this work is publicly available at
\url{https://github.com/PV-314/hypgeom/any-b}. The author is very
happy to help interested readers who have any questions, problems or
suggestions for the use of this code.

\section{Diophantine Approximation via Hypergeometric Functions}
\label{sect:dio}

In this section, we collect notation and statements from Section~2 of \cite{V5}
that we will need in our proof here.

The following lemma is Lemma~\lemDio{} from \cite{V5}. It is a variation of the
so-called ``folklore lemma'' that can be used to obtain irrationality measures
to a number $\theta$ from a sequence of good rational approximations to $\theta$.

\begin{lemma}
\label{lem:2.1}
Let $\theta \in \bbC$ and let $\bbK$ be an imaginary quadratic field. Suppose that
there exist $k_{0},\ell_{0} > 0$ and $E,Q > 1$ such that for all non-negative
integers $r$, there are algebraic integers $p_{r}$ and $q_{r}$ in $\bbK$ with
$\left| q_{r} \right| < k_{0}Q^{r}$ and
$\left| q_{r} \theta - p_{r} \right| \leq \ell_{0}E^{-r}$ satisfying 
$p_{r}q_{r+1} \neq p_{r+1}q_{r}$.

For any algebraic integers $p$ and $q$ in $\bbK$, let $r_{0}$ be the smallest
positive integer such that $\left( Q-1/E  \right)\ell_{0}|q|/\left( Q-1 \right)<cE^{r_{0}}$,
where $0<c<1$.

\noindent
{\rm (a)} We have 
\[
\left| q\theta - p \right| 
> \frac{1-c/E}{k_{0}Q^{r_{0}+1}}.
\]

\noindent
{\rm (b)} When $p/q \neq p_{r_{0}}/q_{r_{0}}$, we have
\[
\left| q\theta - p \right| 
> \frac{1-c}{k_{0}Q^{r_{0}}}.
\]
\end{lemma}

\subsection{Construction of Approximations}
\label{subsect:const}

Let $t'$, $u_{1}$ and $u_{2}$ be rational integers with $t'<0$ such that
$u=\left( u_{1}+u_{2}\sqrt{t'} \right)/2$ be an algebraic integer in
$\bbK=\bbQ \left( \sqrt{t'} \right)$ with $\sigma(u)=\left( u_{1}-u_{2}\sqrt{t'} \right)/2$
as its algebraic (and complex) conjugate. Put $\omega = u/\sigma(u)$ and write
$\omega=e^{i\varphi}$, where $-\pi<\varphi \leq \pi$. For any real number
$\nu$, we shall put $\omega^{\nu}= e^{i\nu\varphi}$ -- unless otherwise stated,
we will use this convention throughout this paper.

Suppose that $\alpha$, $\beta$ and $\gamma$ are complex numbers and $\gamma$ is
not a non-positive integer. We denote by ${}_{2}F_{1}(\alpha, \beta; \gamma; z)$
the classical (or Gauss) hypergeometric function of the complex variable $z$.
For integers $m$ and $n$ with $0 < m < n$, $(m,n) = 1$ and $r$ a non-negative
integer, put $\nu=m/n$ and
\[
X_{m,n,r}(z)={}_{2}F_{1}(-r-\nu, -r; 1-\nu; z), \quad
Y_{m,n,r}=z^{r}X_{m,n,r} \left(z^{-1} \right)
\]
and
\[
R_{m,n,r}(z)
= (z-1)^{2r+1} \frac{\nu \cdots (r+\nu)}{(r+1) \cdots (2r+1)} 
		   {} _{2}F_{1} \left( r+1-\nu, r+1; 2r+2; 1-z \right).
\]

We let $D_{n,r}$ denote the smallest positive integer such that $D_{n,r} X_{m,n,r}(x) \in \bbZ[x]$
for all $m$ as above. For $d' \in \bbZ$, we define $N_{d',n,r}$ to be the largest
integer such that $\left( D_{n,r}/ N_{d',n,r} \right)X_{m,n,r}\left( 1-\sqrt{d'}\,x \right)
\in \bbZ \left[ \sqrt{d'} \right] [x]$, again for all $m$ as above. We will use
$v_{p}(x)$ to denote the largest power of a prime $p$ which divides
the rational number $x$. We put
\begin{equation}
\label{eq:ndn-defn}
\cN_{d',n} =\prod_{p|n} p^{\min(v_{p}(d')/2, v_{p}(n)+1/(p-1))}.
\end{equation}

In what follows, we shall restrict our attention to $m=1$ and $n=4$, so $\nu=1/4$.

As in equation~(\eqGDefn{}) in \cite{V5}, put
\begin{align}
\label{eq:g-defn}
g_{1}  & = \gcd \left( u_{1}, u_{2} \right), \\
g_{2}  & = \gcd \left( u_{1}/g_{1}, t' \right), \nonumber \\
g_{3}  & = \left\{
			 \begin{array}{ll}
	             1 & \text{if $t' \equiv 1 \pmod{4}$ and $\left( u_{1}-u_{2} \right)/g_{1} \equiv 0 \pmod{2}$}, \\
	             2 & \text{if $t' \equiv 3 \pmod{4}$ and $\left( u_{1}-u_{2} \right)/g_{1} \equiv 0 \pmod{2}$},\\
	             4 & \text{otherwise,}
             \end{array}
             \right. \nonumber \\
g      & = g_{1}\sqrt{g_{2}/g_{3}}. \nonumber
\end{align}

Then we can put
\begin{align}
\label{eq:7}
p_{r} &= \frac{D_{4,r}}{N_{d',4,r}} \left( \frac{u_{1}-u_{2}\sqrt{t'}}{2g} \right)^{r} X_{1,4,r}(\omega),\\
q_{r} &= \frac{D_{4,r}}{N_{d',4,r}} \left( \frac{u_{1}-u_{2}\sqrt{t'}}{2g} \right)^{r} Y_{1,4,r}(\omega)
\hspace*{3.0mm} \text{ and} \nonumber \\
R_{r} &= \frac{D_{4,r}}{N_{d',4,r}} \left( \frac{u_{1}-u_{2}\sqrt{t'}}{2g} \right)^{r} R_{1,4,r}(\omega),
\nonumber
\end{align}
where
\begin{equation}
\label{eq:d-defn}
d'=\left( u-\sigma(u) \right)^{2}/g^{2} = u_{2}^{2}t'/g^{2}.
\end{equation}

As in Subsection~\subSectConstr{}\, of \cite{V5}, $p_{r}$ and $q_{r}$ are algebraic
integers in $\bbQ \left( \sqrt{t'} \right)$ and
\[
q_{r}\omega^{1/4} - p_{r}=R_{r}.
\]

As in equations~(\eqQDefn{}), (\eqKUB{}), (\eqEDefn{}) in \cite{V5} and the
expression for $\ell_{0}$ immediately after them there, we have
\begin{align}
\label{eq:q-defn}
Q &= \frac{\cD_{4} \left| \left| u_{1} \right| + \sqrt{u_{1}^{2}-t'u_{2}^{2}} \right|}{|g|\cN_{d',4}},\\
\label{eq:k-UB}
k_{0} &<0.89, \\
\label{eq:e-defn}
E &= \frac{|g|\cN_{d',4} \left|  \left| u_{1} \right|  + \sqrt{u_{1}^{2}-t'u_{2}^{2}} \right|}{\cD_{4}u_{2}^{2}|t'|} \text{ and} \\
\label{eq:ell0-defn}
\ell_{0} &= \cC_{4,2}|\varphi|=0.2|\varphi|.
\end{align}

\section{Lemmas about $\left( x_{k} \right)_{k=-\infty}^{\infty}$ and $\left( y_{k} \right)_{k=-\infty}^{\infty}$}
\label{sect:prelim}

We start by stating Proposition~\propRep{}(b) of \cite{V5} when $\epsilon$ is a
unit in $\bbQ \left( \sqrt{d} \right)$. Since $N_{\varepsilon}=\pm 1$, we have
$\rad \left( \gcd \left( uN_{\alpha}, N_{\varepsilon} \right) \right)=1$, so the
condition
$f | \left( b^{2} \rad \left( \gcd \left( uN_{\alpha}, N_{\varepsilon} \right) \right) \right)$
there becomes $f|b^{2}$ here. This is the result that permits us to
use the hypergeometric method.

\begin{lemma}
\label{lem:quad-rep}
Let $a \neq 0$, $b>0$ and $d$ be rational integers such that $d$ is not a square.
Put $\alpha=a+b^{2} \sqrt{d}$ and denote $N_{\bbQ(\sqrt{d})/\bbQ}(\alpha)$ by
$N_{\alpha}$. Suppose that $-N_{\alpha}$ is a non-zero square, $x \neq 0$ and $y>0$
are rational integers with
\begin{equation}
\label{eq:quad-rep-assumption}
x+y^{2} \sqrt{d} = \alpha \epsilon^{2},
\end{equation}
where $\epsilon=\left( t+u\sqrt{d} \right)/2 \in \cO_{\bbQ \left( \sqrt{d} \right)}$
with $t$ and $u$ non-zero rational integers, norm $N_{\epsilon}=\pm 1$.

We can write
\begin{align}
\pm f^{2} \left( x + N_{\epsilon}\sqrt{N_{\alpha}} \right)
& = \left( a+\sqrt{N_{\alpha}} \right) \left( r + s\sqrt{\core \left( N_{\alpha} \right)} \right)^{4} 
\quad \text{and} \nonumber \\
\label{eq:fy-rel}
fy &= b \left( r^{2}-\core \left( N_{\alpha} \right)s^{2} \right),
\end{align}
for some integers $f$, $r$ and $s$ satisfying $f \neq 0$ and $f | b^{2}$.
\end{lemma}

We will also need lower bounds for the elements in our sequences.
Recall that $K$ is the largest negative integer such that $y_{K}>b^{2}$.

\begin{lemma}
\label{lem:Y-LB1}
Let the $y_{k}$'s be defined by $\eqref{eq:yk-defn}$ with the notation and assumptions
there. Suppose that $N_{\alpha}<0$.
 
\noindent
{\rm (a)}
For all $k$, $2y_{k}$ is a positive integer.
The sequences $\left( y_{k} \right)_{k \geq 0}$ and 
$\left( y_{K+1}, y_{K}, y_{K-1}, y_{K-2}, \ldots \right)$ are increasing sequences of positive numbers.

\noindent
{\rm (b)} We have
\begin{equation}
\label{eq:yLB3-gen-MOOSE}
y_{k} \geq
\left\{
\begin{array}{ll}
\left( \left| N_{\alpha} \right|u^{2} / \left( 4b^{2} \right) \right) \left( du^{2}-3 \right)^{k-1} & \text{for $k>0$,} \\
\left( \left| N_{\alpha} \right|u^{2} / \left( 4b^{2} \right) \right) \left( du^{2}-3 \right)^{\max(0,K-k)} & \text{for $k<0$.}
\end{array}
\right.
\end{equation}
%
\noindent
{\rm (c)} For $N_{\varepsilon}=1$ or $du^{2} \geq 300$, we have
\begin{equation}
\label{eq:yLB4-gen}
y_{k} \geq
\left\{
\begin{array}{ll}
\left( \left| N_{\alpha} \right|u^{2} / \left( 4b^{2} \right) \right) \left( 0.99du^{2} \right)^{k-1} & \text{for $k>0$,} \\
\left( \left| N_{\alpha} \right|u^{2} / \left( 4b^{2} \right) \right) \left( 0.99du^{2} \right)^{\max(0,K-k)} & \text{for $k<0$.}
\end{array}
\right.
\end{equation}

For any $d>1$, we can replace $0.99$ by $0.4$.
\end{lemma}

\begin{proof}
(a) This is part~(b) of Lemma~\lemLB{} of \cite{V5}.

\vspace*{1.0mm}

(b) This is a sharpening of part~(c) of Lemma~\lemLB{} of \cite{V5}.

The proof follows that of part~(c) of Lemma~\lemLB{} of \cite{V5} up to the
following inequality for $k \geq 0$ near the top of page~311 of \cite{V5}:
\[
y_{k+1}=\left( du^{2} + 2N_{\varepsilon} \right)y_{k}-y_{k-1}
\geq \left( du^{2} + 2N_{\varepsilon}-1 \right)y_{k}.
\]

As there, we get a slightly stronger result if $N_{\varepsilon}=1$.
For $N_{\varepsilon}=-1$, the result is immediate. Similarly, for $k<0$.

\vspace*{1.0mm}

(c) The last statement is part~(c) of Lemma~\lemLB{} of \cite{V5}.

Equation~\eqref{eq:yLB4-gen} follows immediately from part~(b) since for
$du^{2} \geq 300$, we have
$\left( du^{2}-3 \right)/ \left( du^{2} \right) \geq 0.99$.
\end{proof}

We now state a gap principle separating distinct squares in the sequence of
$y_{k}$'s defined by \eqref{eq:yk-defn}.

\begin{lemma}
\label{lem:gap}
Let the $y_{k}$'s be defined by $\eqref{eq:yk-defn}$ with $-N_{\alpha}$
a non-zero square. If $y_{i}$ and $y_{j}$ are
distinct squares with $i,j \neq 0$ and
$y_{j}>y_{i} \geq \max \left( 4\sqrt{\left| N_{\alpha} \right|/d}, b^{2}\left| N_{\alpha} \right|/d \right)$,
then
\[
y_{j} > \frac{57.32d^{2}}{b^{4}\left| N_{\alpha} \right|^{2}} y_{i}^{3}.
\]
\end{lemma}

\begin{proof}
This is Lemma~\lemGap{}(a) of \cite{V5}.
\end{proof}

Now we collect some results that we will need for bounding quantities that
arise in the proof of Theorem~\ref{thm:1.2-seq-new} in the next section
(Section~\ref{sect:proof-thm12}).


\begin{lemma}
\label{lem:zeta}
Let the $y_{k}$'s be defined as in $\eqref{eq:yk-defn}$ with the notation and assumptions
there. Suppose that $N_{\alpha}<0$.

Let $y_{k} \geq 4\sqrt{\left| N_{\alpha} \right|/d}$ be a square and put
\[
\omega_{k}=\left( x_{k} + N_{\varepsilon^{k}}\sqrt{N_{\alpha}} \right)
/\left( x_{k} - N_{\varepsilon^{k}}\sqrt{N_{\alpha}} \right)
=e^{i\varphi_{k}}
\]
with $-\pi < \varphi_{k} \leq \pi$.
Then
\[
\left| \varphi_{k} \right| <\frac{2.29\sqrt{\left| N_{\alpha} \right|}}{\left| x_{k} \right|}
<0.6.
\]
\end{lemma}

\begin{proof}
This is Lemma~\lemTheta{}(a) of \cite{V5}.
\end{proof}

We require a result like Lemma~\lemGNd{} of \cite{V5} for our proof. However,
Lemma~\lemGNd{} of \cite{V5} only applies when $\gcd(a,b)=1$, so we prove the
following weaker, but sufficient, inequality here, which holds more
generally without such a $\gcd(a,b)$ condition.

We also note here that Lemma~\lemGcd{} in \cite{V5} requires the additional hypothesis
that $\gcd (a, b)$ is odd. The proof given in \cite{V5} is valid under this
condition. Without this condition, $\gcd \left( x_{k}, y_{k} \right)
/ \gcd (a,b)=1/2$ can also occur.

\begin{lemma}
\label{lem:gNd4}
Let the sequences $\left( x_{k} \right)_{k=-\infty}^{\infty}$ and
$\left( y_{k} \right)_{k=-\infty}^{\infty}$ be as defined in \eqref{eq:yk-defn},
with the notation and assumptions there. Suppose that $k \neq 0$,
$N_{\alpha}<0$ and that $x_{k}$ and $y_{k}$ are both integers.
Using the notation of Subsection~$\ref{subsect:const}$ with $t'=\core \left( N_{\alpha} \right)$,
$u_{1}=2x_{k}$ and $u_{2}=\pm 2\sqrt{N_{\alpha}/\core\left( N_{\alpha} \right)}$,
we also let $d'=u_{2}^{2}t'/g^{2}$, as defined in
$\eqref{eq:d-defn}$, and $\cN_{d',4}$ be as defined in $\eqref{eq:ndn-defn}$.
Then
\begin{equation}
\label{eq:gn-values}
2 \leq |g|\cN_{d',4} \leq 2\sqrt{\left| N_{\alpha} \right|}.
\end{equation}
\end{lemma}

\begin{proof}
From the notation in the statement of the lemma, we have
\[
d'=u_{2}^{2}t'/g^{2}=\left( 4N_{\alpha}/\core\left( N_{\alpha} \right) \right)\core \left( N_{\alpha} \right)/g^{2}
=4N_{\alpha}/g^{2}.
\]

Hence $\cN_{d',4}=2^{\min \left( v_{2}\left( 4N_{\alpha}/g^{2} \right)/2, 3 \right)}$ and so
$\cN_{d',4}|8$.

Since $g_{1}$ and $g_{2}$ are both integers and $g_{3}|4$, we see that $v_{p} \left( g \right) \geq 0$
for odd primes, $p$. From this and since $p \nmid \cN_{d',4}$ for odd primes, $p$,
we see that $v_{p} \left( |g| \cN_{d',4} \right) \geq 0$ for odd primes, $p$.

For $p=2$, we proceed similarly.

If $\min \left( v_{2}\left( 4N_{\alpha}/g^{2} \right)/2, 3 \right)=
v_{2}\left( 4N_{\alpha}/g^{2} \right)/2=v_{2}\left( 4N_{\alpha} \right)/2-v_{2}(g)$,
then $v_{2} \left( |g|\cN_{d',4} \right)=v_{2} \left( 4N_{\alpha} \right)/2 \geq 1$.

If $\min \left( v_{2}\left( 4N_{\alpha}/g^{2} \right)/2, 3 \right)=3$, then from
$v_{2}(g) \geq -v_{2} \left( g_{3} \right)/2 \geq -1$, we have
$v_{2} \left( |g|\cN_{d',4} \right) \geq 2$.

The lower bound for $|g|\cN_{d',4}$ follows.

For the upper bound, we also consider the two possibilities for $\cN_{d',4}$ separately.

If $v_{2}\left( 4N_{\alpha}/g^{2} \right)/2 \leq 3$, then
\[
v_{2} \left( g \cN_{d',4} \right)= v_{2}(g)+v_{2}\left( 4N_{\alpha}/g^{2} \right)/2
=v_{2}\left( 4N_{\alpha} \right)/2.
\]

If $v_{2}\left( 4N_{\alpha}/g^{2} \right)/2 \geq 4$, then
\[
v_{2} \left( g \cN_{d',4} \right)= v_{2}(g)+3
\leq v_{2}(g) + v_{2}\left( 4N_{\alpha}/g^{2} \right)/2 - 1
=v_{2}\left( 4N_{\alpha} \right)/2 - 1.
\]

Since $g_{1}^{2} | \left( 4 N_{\alpha}/\core \left( N_{\alpha} \right) \right)$
and $g_{2} | \core \left( N_{\alpha} \right)$, it follows that the odd part of
$g^{2}$ (and hence the odd part of $g^{2}\cN_{d',4}^{2}$) is at most $\left| N_{\alpha} \right|$.
The upper bound for $|g|\cN_{d',4}$ now follows.
\end{proof}

\begin{lemma}
\label{lem:admissible-b}
Let $a$, $b$ and $d$ be positive integers with $\gcd \left( a, b^{2} \right)$
squarefree. If $db^{4}-a^{2}=n^{2}$ for some positive integer, $n$, then all of
the prime divisors of $b$ are congruent to $1 \bmod{4}$.
\end{lemma}

\begin{proof}
We know that the prime factorisation of $a^{2}+n^{2}$ is
$2^{k} \prod_{i} p_{i}^{\ell_{i}} \prod_{j} q_{j}^{2m_{j}}$, where $k$ is a
non-negative integer, the $\ell_{i}$'s and $m_{j}$'s are positive integers,
the $p_{i}$'s are primes congruent to $1 \bmod{4}$ and
the $q_{j}$'s are primes congruent to $3 \bmod{4}$. Furthermore, we know that
the $q_{j}^{m_{j}}$'s must divide both $a$ and $n$.

By our assumption that $\gcd \left( a, b^{2} \right)$ is squarefree, we know that
$q_{j} \nmid b$. Otherwise, $q_{j}^{4}| b^{4}$ and hence $q_{j}^{4}| \left( a^{2}+n^{2} \right)$.
Thus $m_{j} \geq 2$ and
$q_{j}^{2}| a$ contradicting $\gcd \left( a, b^{2} \right)$ being squarefree.

We now consider the prime $2$. We need not consider $k<4$, as in these cases
$2^{k}|d$ and $2 \nmid b$. If $k \geq 4$, then both $a$ and $n$
must be divisible by $4$. Otherwise, since the squares modulo $16$ are $0,1,4$
and $9$, we have $a^{2}+n^{2} \equiv 1,2,4,5,8,9,10,13 \pmod{16}$.
If we also have $2|b$, then $4|\gcd \left( a, b^{2} \right)$,
contradicting our assumption. Hence $2 \nmid b$. Hence all of the prime divisors
of $b$ are congruent to $1 \bmod{4}$.
\end{proof}

\section{Proof of Theorem~\ref{thm:1.2-seq-new}}
\label{sect:proof-thm12}

Following the same technique as in \cite{V5,V6}, we break the proof into several
parts. There are several algebraic manipulations
in what follows. To check these up to and including Subsection~\ref{subsect:thm12-step-iv},
we wrote some Maple code. It can be found in the file \verb!bAll-nASqr-proof-calcs.txt!
in the \verb!maple! subdirectory of the github url provided at the end of
Section~\ref{sect:intro}.

\subsection{Prerequisites}
\label{subsect:preq}

To prove Theorem~\ref{thm:1.2-seq-new}, we assume there are three distinct squares,
\begin{equation}
\label{eq:ass-7aa}
y_{k_{3}}>y_{k_{2}}>y_{k_{1}} \geq \max \left( 4\sqrt{\left| N_{\alpha} \right|/d}, b^{2}\left| N_{\alpha} \right|/d \right)
\end{equation}
satisfying
\begin{equation}
\label{eq:ass-5aa-xxx}
k_{1} \geq 1 \text{ or } k_{1} \leq K,
\end{equation}
recalling that $K$ is the largest negative integer such that $y_{K} > b^{2}$ and
obtain a contradiction for $y_{k_{1}}$ sufficiently large. Also recall that by
our choice of $\alpha$, there are no squares among the $y_{k}$ with $K+1 \leq k \leq 0$
distinct from $y_{0}$. The condition in
\eqref{eq:ass-7aa} arises from our Gap Principle Lemma above (Lemma~\ref{lem:gap}).
The condition in \eqref{eq:ass-5aa-xxx} ensures that the conditions to apply
both our Gap Principle lemma and Lemma~\ref{lem:quad-rep} are met.

In what follows (e.g., for \eqref{eq:Q-LB1} below), we will also require that
\begin{equation}
\label{eq:ass-3aa}
y_{k_{1}} \geq \max \left( 4\left| N_{\alpha} \right|/\sqrt{d}, |g|\cN_{d',4}/\sqrt{d} \right).
\end{equation}

From the upper bound in \eqref{eq:gn-values} in Lemma~\ref{lem:gNd4}, we have
$|g|\cN_{d',4} \leq 2\sqrt{\left| N_{\alpha} \right|}<4 \left| N_{\alpha} \right|$.
Hence \eqref{eq:ass-3aa} holds if $y_{k_{1}} \geq 4\left| N_{\alpha} \right|/\sqrt{d}$.
Combining this with \eqref{eq:ass-7aa}, we find that \eqref{eq:ass-7aa} and
\eqref{eq:ass-3aa} both hold if
\begin{equation}
\label{eq:C-2}
y_{k_{3}}>y_{k_{2}}>y_{k_{1}} \geq \max \left( 4\left| N_{\alpha} \right|/\sqrt{d}, b^{2}\left| N_{\alpha} \right|/d \right).
\end{equation}

As in Lemma~\ref{lem:zeta}, we put
$\omega_{k_{1}}=\left( x_{k_{1}}+N_{\varepsilon^{k_{1}}}\sqrt{N_{\alpha}} \right)
/\left( x_{k_{1}}-N_{\varepsilon^{k_{1}}}\sqrt{N_{\alpha}} \right)$
and let $\zeta_{4}$ be the $4$-th root of unity such that
\[
\left| \omega_{k_{1}}^{1/4} - \zeta_{4} \frac{x-y\sqrt{\core \left( N_{\alpha} \right)}}{x+y\sqrt{\core \left( N_{\alpha} \right)}} \right|
= \min_{0 \leq j \leq 3} \left( \left| \omega_{k_{1}}^{1/4} - e^{2j\pi i/4} \frac{x-y\sqrt{\core \left( N_{\alpha} \right)}}{x+y\sqrt{\core \left( N_{\alpha} \right)}} \right| \right),
\]
where $x+y\sqrt{\core \left( N_{\alpha} \right)}=\left(  r_{k_{1}}-s_{k_{1}}\sqrt{\core \left( N_{\alpha} \right)} \right)
\left( r_{k_{3}}+s_{k_{3}}\sqrt{\core \left( N_{\alpha} \right)} \right)$ with
$\left( r_{k_{1}}, s_{k_{1}} \right)$ and $\left( r_{k_{3}}, s_{k_{3}} \right)$
as in Lemma~\ref{lem:quad-rep}, which are associated with $\left( x_{k_{1}}, y_{k_{1}} \right)$
and $\left( x_{k_{3}}, y_{k_{3}} \right)$, respectively.
We can take $\zeta_{4} \in \bbQ \left( \sqrt{\core \left( N_{\alpha} \right)} \right)$.
This is immediate here since $-N_{\alpha}$ is a perfect square.
This is important for us here as $\zeta_{4} \left( x-y\sqrt{\core \left( N_{\alpha} \right)} \right)
/ \left( x+y\sqrt{\core \left( N_{\alpha} \right)} \right)$
must be in an imaginary quadratic field in order to apply Lemma~\ref{lem:2.1}
to obtain a lower bound for the rightmost quantity in \eqref{eq:29} below.

From (\eqOmegaUB{}) in the proof of Lemma~\lemTheta{}(b) in \cite{V5}, we have
\[
\left| \omega_{k_{1}}^{1/4} - \zeta_{4} \frac{x-y\sqrt{\core\left( N_{\alpha} \right)}}{x+y\sqrt{\core\left( N_{\alpha} \right)}} \right|
<0.127.
\]
This is where we need \eqref{eq:ass-7aa} (and hence \eqref{eq:C-2}) and
\eqref{eq:ass-5aa-xxx} above.

Thus we can apply Lemma~\lemUB{}(a) in \cite{V5} with $c_{1}=0.127$ to find that
\begin{equation}
\label{eq:29}
\frac{2\sqrt{\left| N_{\alpha} \right|}}{\sqrt{d} \, y_{k_{3}}}
= \left| \omega_{k_{1}} - \left( \frac{x-y\sqrt{\core \left( N_{\alpha} \right)}}{x+y\sqrt{\core \left( N_{\alpha} \right)}} \right)^{4} \right|
> 3.959 \left| \omega_{k_{1}}^{1/4} - \zeta_{4} \frac{x-y\sqrt{\core \left( N_{\alpha} \right)}}{x+y\sqrt{\core \left( N_{\alpha} \right)}} \right|.
\end{equation}

The equality on the left-hand side is from the equalities in (\eqOmegaMinusStuffEquality{})
in the proof of Lemma~\lemTheta{}(b) in \cite{V5}.


We need to derive a lower bound for the far-right quantity in \eqref{eq:29}.
To do so, we shall use the lower bounds in Lemma~\ref{lem:2.1} with
a sequence of good approximations $p_{r}/q_{r}$ obtained from the hypergeometric
functions. So we collect here the required quantities.

Since $y_{k_{1}} \geq 4\sqrt{\left| N_{\alpha} \right|/d}$ (from \eqref{eq:ass-7aa}),
we obtain
\begin{equation}
\label{eq:25}
x_{k_{1}}^{2} = dy_{k_{1}}^{2}+N_{\alpha}
=dy_{k_{1}}^{2} \left( 1 + \frac{N_{\alpha}}{dy_{k_{1}}^{2}} \right)
\geq 0.9375dy_{k_{1}}^{2}.
\end{equation}

So
\begin{equation}
\label{eq:26}
\sqrt{x_{k_{1}}^{2}-N_{\alpha}}= \sqrt{dy_{k_{1}}^{2}}
< 1.04 x_{k_{1}}.
\end{equation}

Using the notation of Subsection~\ref{subsect:const}, let $t'=\core \left( N_{\alpha} \right)$,
$u_{1}=2x_{k_{1}}$, $u_{2}=2\sqrt{N_{\alpha}/\core \left( N_{\alpha} \right)}$ and
$d'$ is as defined in \eqref{eq:d-defn}.

Recall equation~(\eqELB{}) from \cite{V5} with $k$ there being $k_{1}$ here:
\begin{equation}
\label{eq:E-LB1}
E > \frac{0.1832|g|\cN_{d',4}\sqrt{d} \, y_{k_{1}}}{\left| N_{\alpha} \right|}.
\end{equation}

From $|g|\cN_{d',4} \geq 2$ and $y_{k_{1}} \geq 4\left| N_{\alpha} \right|/\sqrt{d}$
we have $E>1$, as required for its use with Lemma~\ref{lem:2.1}.

Similarly, using \eqref{eq:25} and $y_{k_{1}} \geq |g|\cN_{d',4}/\sqrt{d}$, we have
\begin{equation}
\label{eq:Q-LB1}
Q > \frac{2e^{1.68}\left( 1+\sqrt{0.9375} \right) \sqrt{d}y_{k_{1}}}{|g|\cN_{d',4}}
>\frac{21.12\sqrt{d} \, y_{k_{1}}}{|g|\cN_{d',4}}
\geq 21.12.
\end{equation}
For \eqref{eq:E-LB1} and \eqref{eq:Q-LB1}, we have used \eqref{eq:ass-3aa} (and
hence \eqref{eq:C-2}).

From equation~(\eqQUB{}) in \cite{V5}, we have
\begin{equation}
\label{eq:Q-UB2}
Q <\frac{21.47\sqrt{d} \, y_{k_{1}}}{|g|\cN_{d',4}}.
\end{equation}

Writing $\omega_{k_{1}}=e^{i\varphi_{k_{1}}}$, with $-\pi < \varphi_{k_{1}} \leq \pi$,
from \eqref{eq:ell0-defn} and Lemma~\ref{lem:zeta}(a), we can take
\begin{equation}
\label{eq:ell-UB}
\ell_{0}=0.2 \left| \varphi_{k_{1}} \right|<0.46\sqrt{\left| N_{\alpha} \right|}/\left| x_{k_{1}} \right|.
\end{equation}

Also from Lemma~\ref{lem:zeta}, we have $\left| \varphi_{k_{1}} \right|<0.6$,
so the condition $\left| \omega_{k_{1}}-1 \right|<1$ required in Lemma~\lemHypg{}\, of \cite{V5} to
apply the hypergeometric method is satisfied.

Let $q=x+y\sqrt{\core \left( N_{\alpha} \right)}=\left( r_{k_{1}}-s_{k_{1}}\sqrt{\core \left( N_{\alpha} \right)} \right)
\left( r_{k_{3}}+s_{k_{3}}\sqrt{\core \left( N_{\alpha} \right)} \right)$
and $p=x-y\sqrt{\core \left( N_{\alpha} \right)}$. Recall from (\eqQAbs{}) in
\cite{V5} that
\begin{equation}
\label{eq:qAbs}
|q| = \frac{\sqrt{f_{k_{1}}f_{k_{3}}} \left( y_{k_{1}}y_{k_{3}} \right)^{1/4}}{b}.
\end{equation}

From Lemma~\ref{lem:quad-rep}, we have
\begin{equation}
\label{eq:f1f3UB}
f_{k_{1}}f_{k_{3}} \leq b^{4}.
\end{equation}

We are now ready to deduce the required contradiction from the assumptions above.
With $r_{0}$ as in Lemma~\ref{lem:2.1}, we follow the same process as in \cite{V5}.
In each of the four main steps, we will obtain an upper bound for $y_{k_{1}}$, which
must hold in order for the conditions defining that step to hold. Thus we get the
desired lower bound for $y_{k_{1}}$ in the statement of Theorem~\ref{thm:1.2-seq-new}
(since we want none of these four steps to be possible).

\subsection{$r_{0}=1$ and $\zeta_{4}p/q \neq p_{1}/q_{1}$ for all $4$-th roots of unity, $\zeta_{4}$}
\label{subsect:thm12-step-i}

We start by determining an upper bound for $y_{k_{3}}$ for all $r_{0} \geq 1$ when
$\zeta_{4}p/q \neq p_{r_{0}}/q_{r_{0}}$, since we will also need such a result
in Subsection~\ref{subsect:thm12-step-iii}.

From the equality in \eqref{eq:29}, along with Lemma~\ref{lem:2.1}(b) and \eqref{eq:qAbs},
we have
\begin{equation}
\label{eq:32}
\frac{2\sqrt{\left| N_{\alpha} \right|}}{\sqrt{d} \, y_{k_{3}}}
> 3.959 \left| \omega_{k_{1}}^{1/4} - \zeta_{4} \frac{x-y\sqrt{\core \left( N_{\alpha} \right)}}{x+y\sqrt{\core \left( N_{\alpha} \right)}} \right|
> \frac{3.959(1-c)b}{k_{0}Q^{r_{0}}\sqrt{f_{k_{1}}f_{k_{3}}} \left( y_{k_{1}}y_{k_{3}} \right)^{1/4}}.
\end{equation}

Applying \eqref{eq:k-UB} and \eqref{eq:Q-UB2} to \eqref{eq:32}, we obtain
\[
\frac{2\sqrt{\left| N_{\alpha} \right|}}{\sqrt{d} \, y_{k_{3}}}
> \frac{3.959(1-c)b}
  {0.89 \left( 21.47\sqrt{d} \, y_{k_{1}} / \left( |g|\cN_{d',4} \right) \right)^{r_{0}}\sqrt{f_{k_{1}}f_{k_{3}}} \left( y_{k_{1}}y_{k_{3}} \right)^{1/4}}.
\]

After taking the fourth power of both sides and rearranging, we find that
\begin{equation}
\label{eq:y2UB-step1}
\left( N_{\alpha}f_{k_{1}}f_{k_{3}} \right)^{2}\left( \frac{0.45}{(1-c)b} \right)^{4}
\left( \frac{21.47}{|g|\cN_{d',4}} \right)^{4r_{0}}
d^{2r_{0}-2} y_{k_{1}}^{4r_{0}+1}
>y_{k_{3}}^{3}.
\end{equation}

Specialising to the case when $r_{0}=1$ and using $|g|\cN_{d',4} \geq 2$ from the
lower bound in \eqref{eq:gn-values}, we have
\begin{equation}
\label{eq:y2UB-step1-b}
y_{k_{3}}^{3} < 545(b(1-c))^{-4} \left( N_{\alpha}f_{k_{1}}f_{k_{3}} \right)^{2}y_{k_{1}}^{5}.
\end{equation}

Applying \eqref{eq:f1f3UB}, we obtain
\begin{equation}
\label{eq:y3UB}
y_{k_{3}}^{3} < 545(1-c)^{-4} N_{\alpha}^{2}b^{4}y_{k_{1}}^{5}.
\end{equation}

We will now combine \eqref{eq:y3UB} with the gap principle in Lemma~\ref{lem:gap}
to show that this case cannot occur. Applying Lemma~\ref{lem:gap} twice, we
have
\[
y_{k_{3}}> \frac{57.32d^{2}}{b^{4} \left| N_{\alpha} \right|^{2}} y_{k_{2}}^{3}
> \left( \frac{57.32d^{2}}{b^{4} \left| N_{\alpha} \right|^{2}} \right)^{4} y_{k_{1}}^{9}.
\]

Combining the lower bound that this provides for $y_{k_{3}}^{3}$ with the upper
bound for $y_{k_{3}}^{3}$ in \eqref{eq:y3UB} and cancelling the common factor
of $y_{k_{1}}^{5}$ on both sides, we obtain
\[
\left( \frac{57.32d^{2}}{b^{4} \left| N_{\alpha} \right|^{2}} \right)^{12}y_{k_{1}}^{22}
< 545(1-c)^{-4}N_{\alpha}^{2}b^{4},
\]
provided that \eqref{eq:ass-7aa} holds, which holds if \eqref{eq:C-2} holds.

We need \eqref{eq:ass-5aa-xxx} here, since in Lemma~\ref{lem:gap} (which we
use above), we do not permit $k_{1}, k_{3}=0$.

Rearranging, taking $22$-nd roots of both sides, using $c=0.75$ (a choice we
will justify in Subsection~\ref{subsect:thm12-step-iii}) and
$\left( 545(1-0.75)^{-4}/57.32^{12} \right)^{1/22}<0.19$, we find that
\begin{equation}
\label{eq:ykUB-step-i}
\frac{0.19 b^{26/11}\left| N_{\alpha} \right|^{13/11}}{d^{12/11}}>y_{k_{1}}
\end{equation}
must hold in this case.

\subsection{$r_{0}=1$ and $\zeta_{4}p/q = p_{1}/q_{1}$ for some $4$-th root of unity, $\zeta_{4}$}
\label{subsect:thm12-step-ii}

As in Subsection~\ref{subsect:thm12-step-i}, we start by proving an upper bound for
$y_{k_{3}}$ that holds for all $r_{0} \geq 1$ with $\zeta_{4}p/q=p_{r_{0}}/q_{r_{0}}$
for some $4$-th root of unity, $\zeta_{4}$.

From (\eqYLUBUsingYK{}) in \cite{V5}
with $k$ and $\ell$ there set to $k_{1}$ and $k_{3}$, respectively, we have
\begin{equation}
\label{eq:yEll-UB-step2}
1.73r_{0}^{1/2} \left(4\frac{d}{\left| N_{\alpha} \right|} \right)^{r_{0}}y_{k_{1}}^{2r_{0}+1}
>y_{k_{3}}.
\end{equation}

Specialising this inequality to the case of $r_{0}=1$, and assuming \eqref{eq:ass-5aa-xxx}
and \eqref{eq:C-2}, then we can
apply the gap principle in Lemma~\ref{lem:gap} twice to show that
\[
6.92 \frac{d}{\left| N_{\alpha} \right|} y_{k_{1}}^{3}>y_{k_{3}}
>57.32 \frac{d^{2}}{b^{4} \left| N_{\alpha} \right|^{2}} y_{k_{2}}^{3}
> \left( 57.32 \frac{d^{2}}{b^{4} \left| N_{\alpha} \right|^{2}} \right)^{4} y_{k_{1}}^{9}.
\]

Rearranging this, we find that
\[
\frac{10^{-6}\left| N_{\alpha} \right|^{7}b^{16}}{d^{7}}>y_{k_{1}}^{6}
\]
must hold. Taking $6$-th roots,
\begin{equation}
\label{eq:ykUB-step-ii}
\frac{\left| N_{\alpha} \right|^{7/6}b^{8/3}}{10d^{7/6}}
>y_{k_{1}}
\end{equation}
must hold.

\subsection{$r_{0}>1$, $\zeta_{4}p/q \neq p_{r_{0}}/q_{r_{0}}$ for all $4$-th
roots of unity, $\zeta_{4}$}
\label{subsect:thm12-step-iii}

Here we establish a stronger gap principle for $y_{k_{1}}$ and $y_{k_{3}}$ than
the one in Lemma~\ref{lem:gap}. We then use this with the upper bound for $y_{k_{3}}$
in \eqref{eq:y2UB-step1} to obtain a contradiction.

We start by deriving a lower bound for $y_{k_{3}}$ that holds in both this step
and in the next step.

From the definition of $r_{0}$ in Lemma~\ref{lem:2.1},
along with \eqref{eq:Q-LB1} and $E>1$, we have
\begin{equation}
\label{eq:qLB-step3}
|q| \geq \frac{c(Q-1)}{\ell_{0}(Q-1/E)}E^{r_{0}-1}
> 0.952cE^{r_{0}-1}/\ell_{0}.
\end{equation}

Recall that $|q|=\sqrt{f_{k_{1}}f_{k_{3}}} \left( y_{k_{1}}y_{k_{3}} \right)^{1/4}/b$
by \eqref{eq:qAbs}. Thus
\[
\left( y_{k_{1}}y_{k_{3}} \right)^{1/4}
> \frac{0.952bcE^{r_{0}-1}}{\ell_{0}\sqrt{f_{k_{1}}f_{k_{3}}}}.
\]

Applying \eqref{eq:E-LB1} and \eqref{eq:ell-UB}, and then \eqref{eq:25}, to this
inequality, we obtain
\begin{align*}
\left( y_{k_{1}}y_{k_{3}} \right)^{1/4}
&> \frac{0.952bc \left| x_{k_{1}} \right|}{0.46\sqrt{\left| N_{\alpha} \right| f_{k_{1}}f_{k_{3}}}}
\left( \frac{0.1832|g|\cN_{d',4}\sqrt{d} \, y_{k_{1}}}{\left| N_{\alpha} \right|} \right)^{r_{0}-1} \\
&> \frac{2bc\sqrt{d} \, y_{k_{1}}}{\sqrt{\left| N_{\alpha} \right|f_{k_{1}}f_{k_{3}}}}
  \left( \frac{0.1832|g|\cN_{d',4}\sqrt{d} \, y_{k_{1}}}{\left| N_{\alpha} \right|} \right)^{r_{0}-1}.
\end{align*}

Taking the fourth power of both sides and rearranging, we find that
\begin{equation}
\label{eq:y2LB-step3}
y_{k_{3}} > \left( \frac{10.91bc}{|g|\cN_{d',4}} \sqrt{\frac{\left| N_{\alpha} \right|}{f_{k_{1}}f_{k_{3}}}} \right)^{4}
           \left( \frac{0.1832 |g|\cN_{d',4}}{\left| N_{\alpha} \right|} \right)^{4r_{0}} d^{2r_{0}}
y_{k_{1}}^{4r_{0}-1}.
\end{equation}

With this lower bound for $y_{k_{3}}$, we now focus for the rest of this subsection
on when $\zeta_{4}p/q \neq p_{r_{0}}/q_{r_{0}}$ for all $4$-th
roots of unity, $\zeta_{4}$.

We now take the third power of both sides of this inequality and combine it
with the upper bound for $y_{k_{3}}^{3}$ in \eqref{eq:y2UB-step1}, finding that
\begin{align}
\label{eq:y1UB-step3-a}
&
\left( \left| N_{\alpha} \right| f_{k_{1}}f_{k_{3}} \right)^{2} \left( \frac{0.45}{b(1-c)} \right)^{4}
\left( \frac{21.47}{|g|\cN_{d',4}} \right)^{4r_{0}}
d^{2r_{0}-2} y_{k_{1}}^{4r_{0}+1} \\
& >
\left( \frac{10.91bc}{|g|\cN_{d',4}} \sqrt{\frac{\left| N_{\alpha} \right|}{f_{k_{1}}f_{k_{3}}}} \right)^{12}
           \left( \frac{0.1832 |g|\cN_{d',4}}{\left| N_{\alpha} \right|} \right)^{12r_{0}}
		   d^{6r_{0}} y_{k_{1}}^{12r_{0}-3} \nonumber
\end{align}
must hold for us to be in this case.

As in Subsection~4.4 of \cite{V5}, $c^{12}(1-c)^{4}$ is monotonically increasing
for $0<c \leq 0.75$. So we put $c=0.75$ and have $c^{12}(1-c)^{4}>0.000124$.
Applying this to \eqref{eq:y1UB-step3-a} and simplifying, we find that if
\eqref{eq:y1UB-step3-a} holds, then
\[
\left( f_{k_{1}}f_{k_{3}} \right)^{8}
      \frac{0.00143\left| N_{\alpha} \right|^{2} \left( |g|\cN_{d',4} \right)^{4}}{b^{16}d^{4}}
      \left( \frac{1.2194 \cdot 10^{7}}
                  {\left( |g|\cN_{d',4} \right)^{8}} \right)^{2r_{0}-1}
> \left( \frac{y_{k_{1}}^{4}d^{2}}{\left| N_{\alpha} \right|^{6}} \right)^{2r_{0}-1}
\]
must hold.

That is,
\begin{equation}
\label{eq:y1UB-step-iiib}
\left( f_{k_{1}}f_{k_{3}} \right)^{8}
\frac{0.00143\left| N_{\alpha} \right|^{2} \left( |g|\cN_{d',4} \right)^{4}}{b^{16}d^{4}}
> \left( \frac{y_{k_{1}}^{4}d^{2}\left( |g|\cN_{d',4} \right)^{8}}
{59.2^{4} \left| N_{\alpha} \right|^{6}} \right)^{2r_{0}-1}
\end{equation}
must hold.
Applying \eqref{eq:f1f3UB}, the left-hand side of \eqref{eq:y1UB-step-iiib}
is at most
\[
b^{16}\frac{0.00143\left| N_{\alpha} \right|^{2} \left( |g|\cN_{d',4} \right)^{4}}{d^{4}}.
\]

So, since $0.00143<0.195^{4}$, if \eqref{eq:y1UB-step-iiib} holds, then
\begin{equation}
\label{eq:ykUB-step-iiia}
\frac{59.2 \cdot 0.195^{1/(2r_{0}-1)}\left| N_{\alpha} \right|^{(3/2)+1/(2(2r_{0}-1))} b^{4/(2r_{0}-1)}}
{d^{(1/2)+1/(2r_{0}-1)} \left( |g|\cN_{d',4} \right)^{2-1/(2r_{0}-1)}}
>y_{k_{1}}
\end{equation}
must also hold.

We can simplify this upper bound for $y_{k_{1}}$.
From Lemma~\ref{lem:gNd4}, we know that
$|g|\cN_{d',4} \geq 2$, so $\left( |g|\cN_{d',4} \right)^{-2+1/(2r_{0}-1)} \leq
\left( |g|\cN_{d',4} \right)^{1/(2r_{0}-1)}/4$. Applying this and $r_{0} \geq 2$,
we obtain
\begin{align*}
& \frac{59.2 \cdot 0.195^{1/(2r_{0}-1)} b^{4/(2r_{0}-1)} \left| N_{\alpha} \right|^{(3/2)+1/(2(2r_{0}-1))}}
{d^{(1/2)+1/(2r_{0}-1)} \left( |g|\cN_{d',4} \right)^{2-1/(2r_{0}-1)}} \\
<&
\frac{59.2 \cdot 0.39^{1/(2r_{0}-1)} b^{4/(2r_{0}-1)} \left| N_{\alpha} \right|^{(3/2)+1/(2(2r_{0}-1))}}
{4d^{1/2}} \\
<& \frac{14.8b^{4/3} \left| N_{\alpha} \right|^{5/3}}{d^{1/2}},
\end{align*}
where we have also used $d \geq 2$ in the last inequality.

Hence if we are in this case, then
\begin{equation}
\label{eq:ykUB-step-iii}
\frac{14.8b^{4/3} \left| N_{\alpha} \right|^{5/3}}{d^{1/2}}>y_{k_{1}}
\end{equation}
must hold.

\subsection{$r_{0}>1$ and $\zeta_{4}p/q = p_{r_{0}}/q_{r_{0}}$ for some $4$-th root of unity, $\zeta_{4}$}
\label{subsect:thm12-step-iv}

We now combine our upper bound for $y_{k_{3}}$ in \eqref{eq:yEll-UB-step2} with our
lower bound for $y_{k_{3}}$ in \eqref{eq:y2LB-step3}. Thus
\[
1.73r_{0}^{1/2} \left(4\frac{d}{\left| N_{\alpha} \right|} \right)^{r_{0}}y_{k_{1}}^{2r_{0}+1}
> \left( \frac{10.91bc}{|g|\cN_{d',4}} \sqrt{\frac{\left| N_{\alpha} \right|}{f_{k_{1}}f_{k_{3}}}} \right)^{4}
  \left( \frac{0.1832 |g|\cN_{d',4}}{\left| N_{\alpha} \right|} \right)^{4r_{0}} d^{2r_{0}}
  y_{k_{1}}^{4r_{0}-1}
\]
and so
\[
1.73r_{0}^{1/2}
> \left( \frac{10.91bc}{|g|\cN_{d',4}} \sqrt{\frac{\left| N_{\alpha} \right|}{f_{k_{1}}f_{k_{3}}}} \right)^{4}
  \left( \frac{0.1832^{4} |g|^{4}\cN_{d',4}^{4}}{4\left| N_{\alpha} \right|^{3}} \right)^{r_{0}} d^{r_{0}}
  y_{k_{1}}^{2r_{0}-2}.
\]

We showed in Subsection~\subSectStepiv{} of \cite{V5} that $0.1832^{4r_{0}}/r_{0}^{1/2}>0.175^{4r_{0}}$.
Applying this, along with $c=0.75$ and collecting
the terms taken to the power $r_{0}-1$, yields
\begin{equation}
1
> \frac{0.607b^{4}d}{f_{k_{1}}^{2}f_{k_{3}}^{2} \left| N_{\alpha} \right|}
  \left( \frac{0.0002344 |g|^{4}\cN_{d',4}^{4}}{\left| N_{\alpha} \right|^{3}} dy_{k_{1}}^{2} \right)^{r_{0}-1}.
\end{equation}

We now proceed similarly to the way we did in Subsection~\ref{subsect:thm12-step-iii}.

Applying \eqref{eq:f1f3UB}, we have
\[
\frac{1.648b^{4}\left| N_{\alpha} \right|}{d}
> \left( \frac{0.0002344 |g|^{4}\cN_{d',4}^{4}}{\left| N_{\alpha} \right|^{3}} dy_{k_{1}}^{2} \right)^{r_{0}-1}.
\]

In order for this inequality to hold, we must have
\begin{equation}
\label{eq:ykUB-step-iva}
65.32\frac{b^{2/(r_{0}-1)}\left| N_{\alpha} \right|^{(3/2)+1/(2(r_{0}-1))} \cdot 1.648^{1/(2(r_{0}-1))}}
{|g|^{2}\cN_{d',4}^{2} d^{(1/2)+1/(2(r_{0}-1))}}
>y_{k_{1}}.
\end{equation}

We now simplify this upper bound. From Lemma~\ref{lem:gNd4}, we know that
$|g|\cN_{d',4} \geq 2$. Applying this and $r_{0} \geq 2$, we obtain
\begin{align*}
& \frac{65.32 \cdot 1.648^{1/(2(r_{0}-1))} b^{2/(r_{0}-1)}\left| N_{\alpha} \right|^{(3/2)+1/(2(r_{0}-1))}}
{d^{(1/2)+1/(2(r_{0}-1))} \left( |g|\cN_{d',4} \right)^{2}} \\
\leq & \frac{16.33 \cdot 1.648^{1/(2(r_{0}-1))} b^{2/(r_{0}-1)}\left| N_{\alpha} \right|^{(3/2)+1/(2(r_{0}-1))}}
{d^{(1/2)+1/(2(r_{0}-1))}} \\
<& \frac{16.33 b^{2}\left| N_{\alpha} \right|^{2}}{d^{1/2}},
\end{align*}
where we have also used $d \geq 2$ in the last inequality, so $(1.648/d)^{1/(2(r_{0}-1))}<1$.

Hence in this case,
\begin{equation}
\label{eq:ykUB-step-iv}
y_{k_{1}}>\frac{16.33 b^{2}\left| N_{\alpha} \right|^{2}}{d^{1/2}}.
\end{equation}
must hold.

\subsection{Consolidation}
\label{subsect:thm12-consolidation}

Recall that we have assumed that there are three distinct squares,
$y_{k_{3}}>y_{k_{2}}>y_{k_{1}}$ with $k_{1} \geq 1$ or $k_{1} \leq K$.
Combining \eqref{eq:C-2}, \eqref{eq:ykUB-step-i}, \eqref{eq:ykUB-step-ii}, \eqref{eq:ykUB-step-iii}
and \eqref{eq:ykUB-step-iv}, we see that we get a contradiction to our assumption if
\[
y_{k_{1}}> \max \left(
\frac{0.19 b^{26/11}\left| N_{\alpha} \right|^{13/11}}{d^{12/11}},
\frac{b^{8/3}\left| N_{\alpha} \right|^{7/6}}{10d^{7/6}},
\frac{14.8b^{4/3} \left| N_{\alpha} \right|^{5/3}}{d^{1/2}},
\frac{16.33 b^{2}\left| N_{\alpha} \right|^{2}}{d^{1/2}},
\frac{4\left| N_{\alpha} \right|}{\sqrt{d}},
\frac{b^{2}\left| N_{\alpha} \right|}{d}
\right).
\]

Taking the largest constant, the largest exponent on $\left| N_{\alpha} \right|$,
the largest exponent on $b$ and the smallest exponent on $d$, we get a contradiction
if
\[
y_{k_{1}}> \frac{16.33 b^{8/3}\left| N_{\alpha} \right|^{2}}{\sqrt{d}},
\]
proving Theorem~\ref{thm:1.2-seq-new}.

We also present here a more precise version of Theorem~\ref{thm:1.2-seq-new} using
the above lower bound for $y_{k_{1}}$.
This may be of interest to readers. It is also used below in Subsection~\ref{subsect:proof-thm13a}
for the proof of Theorem~\ref{thm:1.3-seq-new}(a).

\begin{proposition}
\label{prop:1.2-seq-new}
Let $a$, $b$ and $d$ be positive integers, where $d$ is not a square, $N_{\alpha}<0$
and $-N_{\alpha}$ is a square. There are at most two distinct squares among the
$y_{k}$'s with $k \geq 1$ or $k \leq K$, and
\[
y_{k} > \max \left(
\frac{0.19 b^{26/11}\left| N_{\alpha} \right|^{13/11}}{d^{12/11}},
\frac{b^{8/3}\left| N_{\alpha} \right|^{7/6}}{10d^{7/6}},
\frac{14.8b^{4/3} \left| N_{\alpha} \right|^{5/3}}{d^{1/2}},
\frac{16.33 b^{2}\left| N_{\alpha} \right|^{2}}{d^{1/2}},
\frac{4\left| N_{\alpha} \right|}{\sqrt{d}},
\frac{b^{2}\left| N_{\alpha} \right|}{d}
\right).
\]
\end{proposition}

\section{Proof of Theorem~\ref{thm:1.3-seq-new}}
\label{sect:proof-thm13}

Here we will combine our bounds on $y_{k}$ in Section~\ref{sect:proof-thm12}
with \eqref{eq:ym1LB} below to obtain Theorem~\ref{thm:1.3-seq-new}. As in
Section~\ref{sect:proof-thm12}, we proceed in separate steps.
Also as in that section (see \eqref{eq:ass-7aa}), here we will assume here that
there are three
distinct squares,
\[
y_{k_{3}}>y_{k_{2}}>y_{k_{1}} \geq \max \left( 4\sqrt{\left| N_{\alpha} \right|/d}, b^{2}\left| N_{\alpha} \right|/d \right)
\]
among the $y_{k}$'s. However here we will assume that
\begin{equation}
\label{eq:ass-5aa}
k_{1} \geq 2 \text{ or } k_{1} \leq K-1.
\end{equation}
In this way, we will obtain a contradiction once $d$ is sufficiently large.

The reason for the stronger condition on $k_{1}$ in \eqref{eq:ass-5aa}, than
that in \eqref{eq:ass-5aa-xxx}, is to
obtain a stronger lower bound for $y_{k_{1}}$ from Lemma~\ref{lem:Y-LB1}(c).

In order to apply Lemma~\ref{lem:Y-LB1}(c), we will also assume that
\begin{equation}
\label{eq:ass-dLB-from-yLB}
d \geq \frac{300}{u^{2}}.
\end{equation}

Applying \eqref{eq:ass-5aa} and \eqref{eq:ass-dLB-from-yLB} to Lemma~\ref{lem:Y-LB1}(c),
we obtain
\begin{equation}
\label{eq:ym1LB}
y_{k_{1}} \geq \frac{0.99|N_{\alpha}|du^{4}}{4b^{2}}.
\end{equation}

By Theorem~1.4 in \cite{V5} and Lemma~\ref{lem:admissible-b}, we see that a
stronger version of Theorem~\ref{thm:1.3-seq-new} holds for $b=1,2,3$ and $4$.
So we may assume that
\begin{equation}
\label{eq:b-LB}
b \geq 5.
\end{equation}


\subsection{$r_{0}=1$ and $\zeta_{4}p/q \neq p_{1}/q_{1}$ for all $4$-th roots of unity, $\zeta_{4}$}
\label{subsect:thm13-step-i}

In Subsection~\ref{subsect:thm12-step-i}, we showed that
\[
y_{k_{1}} < \frac{0.19 b^{26/11}\left| N_{\alpha} \right|^{13/11}}{d^{12/11}}.
\]
This is \eqref{eq:ykUB-step-i} there.
This contradicts \eqref{eq:ym1LB} if
\begin{equation}
\label{eq:ass-6a}
d>\frac{0.882 \left| N_{\alpha} \right|^{2/23} b^{48/23}}{u^{44/23}}.
\end{equation}

It is useful to get a lower bound for $d$ just in terms of $b$ and $u$. Applying
$\left| N_{\alpha} \right| < db^{4}$, we also see that the upper and lower bounds
for $y_{k_{1}}$ in \eqref{eq:ym1LB} and \eqref{eq:ykUB-step-i} contradict each
other if
\begin{equation}
\label{eq:ass-6b}
d>\frac{0.88 b^{8/3}}{u^{44/21}}.
\end{equation}

Hence we cannot be in this step if either \eqref{eq:ass-6a} or \eqref{eq:ass-6b}
holds, along with the earlier assumptions made in Subsection~\ref{subsect:preq}.

\subsection{$r_{0}=1$ and $\zeta_{4}p/q=p_{1}/q_{1}$ for some $4$-th root of unity, $\zeta_{4}$}
\label{subsect:thm13-step-ii}

Combining \eqref{eq:ykUB-step-ii} with \eqref{eq:ym1LB}, we get a contradiction
provided that
\begin{equation}
\label{eq:ass-9a}
d > 0.66\frac{\left| N_{\alpha} \right|^{1/13}b^{28/13}}{u^{24/13}}.
\end{equation}

Applying $\left| N_{\alpha} \right|<db^{4}$, we obtain
\begin{equation}
\label{eq:ass-9b}
d > 0.64\frac{b^{8/3}}{u^{2}}.
\end{equation}

\subsection{$r_{0}>1$, $\zeta_{4}p/q \neq p_{r_{0}}/q_{r_{0}}$ for all $4$-th
roots of unity, $\zeta_{4}$}
\label{subsect:thm13-step-iii}

Combining \eqref{eq:ykUB-step-iiia} and \eqref{eq:ym1LB}, if
\[
\frac{0.99|N_{\alpha}|du^{4}}{4b^{2}}
>\frac{59.2 \cdot 0.195^{1/(2r_{0}-1)}\left| N_{\alpha} \right|^{(3/2)+1/(2(2r_{0}-1))} b^{4/(2r_{0}-1)}}
{d^{(1/2)+1/(2r_{0}-1)} \left( |g|\cN_{d',4} \right)^{2-1/(2r_{0}-1)}},
\]
then the right-hand side of \eqref{eq:y1UB-step-iiib} is greater than \eqref{eq:ykUB-step-iiia}.

Since $|g|\cN_{d',4} \geq 2$, this inequality holds if
\begin{equation}
\label{eq:dLB-step3a}
d >\frac{239.2^{2(2r_{0}-1)/(6r_{0}-1)} \cdot 0.195^{2/(6r_{0}-1)}
\left| N_{\alpha} \right|^{2r_{0}/(6r_{0}-1)} b^{4(2r_{0}+1)/(6r_{0}-1)}}
{u^{8(2r_{0}-1)/(6r_{0}-1)} 2^{2(4r_{0}-3)/(6r_{0}-1)}}.
\end{equation}

For $r_{0} \geq 2$, we have the following:\\
$2r_{0}/ \left( 6r_{0}-1 \right)$ decreases from $4/11=0.3636\ldots$ towards $1/3$,
so $\left| N_{\alpha} \right|^{2r_{0}/(6r_{0}-1)} \leq \left| N_{\alpha} \right|^{4/11}$,\\
$4 \left( 2r_{0}+1 \right)/ \left( 6r_{0}-1 \right)$ decreases from $20/11$
towards $4/3$,
so $b^{4(2r_{0}+1)/(6r_{0}-1)} \leq b^{20/11}$,\\
$8 \left( 2r_{0}-1 \right)/ \left( 6r_{0}-1 \right)$ increases from $24/11$
towards $8/3$,
so $u^{8(2r_{0}-1)/(6r_{0}-1)} \geq u^{24/11}$.

Also, $2 \left( 2r_{0}-1 \right)/ \left( 6r_{0}-1 \right)=2/3-4/\left( 3 \left( 6r_{0}-1 \right) \right)$
and $2 \left( 4r_{0}-3 \right)/ \left( 6r_{0}-1 \right)=4/3-14/\left( 3 \left( 6r_{0}-1 \right) \right)$,
so
\[
\frac{239.2^{2(2r_{0}-1)/(6r_{0}-1)} \cdot 0.195^{2/(6r_{0}-1)}}{2^{2(4r_{0}-3)/(6r_{0}-1)}}
=\frac{239.2^{2/3}}{2^{4/3}}
\left( \frac{0.195^{6} \cdot 2^{14}}{239.2^{4}} \right)^{1/(3(6r_{0}-1))}.
\]

Since $0.195^{6} \cdot 2^{14}/239.2^{4}<1$, we find that
$\left( 0.195^{6} \cdot 2^{14}/239.2^{4} \right)^{1/(3(6r_{0}-1))}<1$. Since
$239.2^{2/3}/2^{4/3}<15.3$, \eqref{eq:dLB-step3a} holds if
\begin{equation}
\label{eq:ass-10a}
d >\frac{15.3\left| N_{\alpha} \right|^{4/11} b^{20/11}}
{u^{24/11}}.
\end{equation}

If we apply $\left| N_{\alpha} \right|<db^{4}$ to \eqref{eq:dLB-step3a}, we obtain
\[
d^{(4r_{0}-1)/(6r_{0}-1)} >\frac{239.2^{2(2r_{0}-1)/(6r_{0}-1)} \cdot 0.195^{2/(6r_{0}-1)} b^{(16r_{0}+4)/(6r_{0}-1)}}
{u^{8(2r_{0}-1)/(6r_{0}-1)} 2^{2(4r_{0}-3)/(6r_{0}-1)}}.
\]

That is,
\[
d>\frac{239.2^{2(2r_{0}-1)/(4r_{0}-1)} \cdot 0.195^{2/(4r_{0}-1)} b^{(16r_{0}+4)/(4r_{0}-1)}}
{u^{8(2r_{0}-1)/(4r_{0}-1)} 2^{2(4r_{0}-3)/(4r_{0}-1)}}.
\]

We have
\begin{align*}
& 239.2^{2(2r_{0}-1)/(4r_{0}-1)} \cdot 0.195^{2/(4r_{0}-1)} / 2^{2(4r_{0}-3)/(4r_{0}-1)}
=239.2^{1-1/(4r_{0}-1)} \cdot 0.195^{2/(4r_{0}-1)} / 2^{2-4/(4r_{0}-1)} \\
= & (239.2/4) \left( 0.195^{2} \cdot 16/239.2 \right)^{1/(4r_{0}-1)}
<60,
\end{align*}
with the maximum value for $r_{0} \geq 2$ occurring as $r_{0} \rightarrow +\infty$,
since $0.195^{2} \cdot 16/239.2<1$. So
\[
d > \frac{60b^{4+8/(4r_{0}-1)}}
{u^{4-4/(4r_{0}-1)}}.
\]

Since $r_{0} \geq 2$, the desired lower bound for $y_{k_{1}}$ holds if
\begin{equation}
\label{eq:ass-10b}
d>\frac{60 b^{36/7}}{u^{24/7}}.
\end{equation}

\subsection{$r_{0}>1$ and $\zeta_{4}p/q = p_{r_{0}}/q_{r_{0}}$ for some $4$-th root of unity, $\zeta_{4}$}
\label{subsect:thm13-step-iv}

We want to show that for $d$ sufficiently large, the upper bound for $y_{k_{1}}$
in \eqref{eq:ykUB-step-iva} is smaller than $0.99\left| N_{\alpha} \right| du^{4}/ \left( 4b^{2} \right)$,
from \eqref{eq:ym1LB}.

So, observing that $65.32/(0.99/4)<264$, we want to show that
\[
d^{(3/2)+1/(2(r_{0}-1))}
>264\frac{b^{2+2/(r_{0}-1)}\left| N_{\alpha} \right|^{(1/2)+1/(2(r_{0}-1))} \cdot 1.648^{1/(2(r_{0}-1))}}
{u^{4}|g|^{2}\cN_{d',4}^{2}}.
\]

This becomes
\begin{align}
\label{eq:dLB1-step-iv}
d> & 264^{(2/3)-2/(3(3r_{0}-2))}\frac{b^{(4/3)+8/(3(3r_{0}-2))}\left| N_{\alpha} \right|^{(1/3)+2/(3(3r_{0}-2))} \cdot 1.648^{1/(3r_{0}-2)}}
{u^{(8/3)-8/(3(3r_{0}-2))}\left( |g|\cN_{d',4} \right)^{(4/3)-4/(3(3r_{0}-2))}} \\
= & \frac{264^{2/3}}{\left( |g|\cN_{d',4} \right)^{4/3}}
\left( \frac{1.648 \left( |g|\cN_{d',4} \right)^{4/3}}{264^{2/3}} \right)^{1/(3r_{0}-2)}
\frac{b^{(4/3)+8/(3(3r_{0}-2))} \left| N_{\alpha} \right|^{(1/3)+2/(3(3r_{0}-2))}}
{u^{(8/3)-8/(3(3r_{0}-2))}}. \nonumber
\end{align}

Since $|g|\cN_{d',4} \geq 2$ and $r_{0} \geq 2$, we see that
$\left( |g|\cN_{d',4} \right)^{(4/3)(1/(3r_{0}-2)-1)}
\leq \left( |g|\cN_{d',4} \right)^{(4/3)(1/4-1)}
= \left( |g|\cN_{d',4} \right)^{-1}$ decreases as $|g|\cN_{d',4}$ increases.
Hence
\[
\frac{264^{2/3}}{\left( |g|\cN_{d',4} \right)^{4/3}}
\left( \frac{1.648 \left( |g|\cN_{d',4} \right)^{4/3}}{264^{2/3}} \right)^{1/(3r_{0}-2)}
\leq \frac{264^{2/3}}{2^{4/3}}
\left( \frac{1.648 \cdot 2^{4/3}}{264^{2/3}} \right)^{1/(3r_{0}-2)}.
\]

We have
$264^{2/3}/2^{4/3}<17$ and
$1.648 \cdot 2^{4/3}/264^{2/3}<1$. We also see that the
biggest exponents on $b$ and $\left| N_{\alpha} \right|$, as well as the smallest
exponent on $u$, all occur for $r_{0}=2$. Hence \eqref{eq:dLB1-step-iv}
holds if
\begin{equation}
\label{eq:dLB2-step-iv}
d>\frac{17b^{2}\left| N_{\alpha} \right|^{1/2}}{u^{2}}.
\end{equation}

We now apply $\left| N_{\alpha} \right|<db^{4}$ to \eqref{eq:dLB1-step-iv}, obtaining
\[
d>\left( \frac{264^{2}}{\left( |g|\cN_{d',4} \right)^{4}} \right)^{(r_{0}-1)/(3r_{0}-2)}
\frac{b^{8r_{0}/(3r_{0}-2)} d^{r_{0}/(3r_{0}-2)} \cdot 1.648^{1/(3r_{0}-2)}}
{u^{8(r_{0}-1)/(3r_{0}-2)}}.
\]

Moving the factor of $d$ to the left-hand side and taking the $\left( 3r_{0}-2 \right)/ \left( 2 \left( r_{0}-1 \right) \right)$-th
root of both sides, this becomes
\[
d>264\frac{b^{4+4/(r_{0}-1)} \cdot 1.648^{1/(2(r_{0}-1))}}
{u^{4}\left( |g|\cN_{d',4} \right)^{2}}.
\]

Since $r_{0} \geq 2$ and $|g|\cN_{d',4} \geq 2$, we obtain
$264 \cdot 1.648^{1/(2(r_{0}-1))}/ \left( |g|\cN_{d',4} \right)^{2}
\leq 264 \cdot 1.648^{1/2}/4<85$, so
\begin{equation}
\label{eq:dLB3-step-iv}
d>85\frac{b^{8}}{u^{4}}.
\end{equation}

If we apply $\left| N_{\alpha} \right| < db^{4}$ to \eqref{eq:dLB2-step-iv} instead
of \eqref{eq:dLB1-step-iv}, then we obtain $290b^{8}/u^{4}$. The improved bound
in \eqref{eq:dLB3-step-iv} will be helpful in the proof of Theorem~\ref{thm:1.3-seq-new}(a).

\subsection{Proof of Theorem~\ref{thm:1.3-seq-new}(b)}
\label{subsect:thm13-step-v}

We now bring together all of the conditions that we have imposed to obtain a
contradiction from the assumption that for $d$ sufficiently large, there are
three distinct squares with $k \geq 2$ or $k \leq K-1$.

First, to be able to apply \eqref{eq:ym1LB}, we need
\begin{equation}
\label{eq:dUB-from-ass-1a}
d \geq \frac{300}{u^{2}}.
\end{equation}

For \eqref{eq:C-2}, we need 
$y_{k_{1}} \geq 4\left| N_{\alpha} \right|/\sqrt{d}$.
By \eqref{eq:ym1LB}, this holds if
\begin{equation}
\label{eq:dUB-from-ass-2a}
d \geq \frac{7b^{4/3}}{u^{8/3}}.
\end{equation}

We also need $y_{k_{1}} \geq b^{2}\left| N_{\alpha} \right|/d$ for \eqref{eq:C-2}.
Using \eqref{eq:ym1LB} again, this holds if
\begin{equation}
\label{eq:dUB-from-ass-3a}
d \geq \frac{2.1b^{2}}{u^{2}}.
\end{equation}

In Subsection~\ref{subsect:thm13-step-i}, we also added the assumption that
$d>0.882 \left| N_{\alpha} \right|^{2/23} b^{48/23}/u^{44/23}$ in \eqref{eq:ass-6a}.

In Subsection~\ref{subsect:thm13-step-ii}, we also added the assumption that
$d>0.66\left| N_{\alpha} \right|^{1/13}b^{28/13}/u^{24/13}$ in \eqref{eq:ass-9a}.

Near the end of Subsection~\ref{subsect:thm13-step-iii}, we imposed the condition that
$d>15.3 \left| N_{\alpha} \right|^{4/11}b^{20/11}/u^{24/11}$ in \eqref{eq:ass-10a}.

Finally, near the end of Subsection~\ref{subsect:thm13-step-iv}, we imposed the condition
that $d>17 \left| N_{\alpha} \right|^{1/2}b^{2}/u^{2}$ in \eqref{eq:dLB2-step-iv}.

It is immediate that \eqref{eq:dUB-from-ass-2a} and \eqref{eq:dUB-from-ass-3a}
follow from this last lower bound for $d$. Similarly, by applying
\eqref{eq:b-LB}, we have $b \geq 5$,
so $17 \left| N_{\alpha} \right|^{1/2}b^{2}/u^{2} \geq 17 \cdot 25/u^{2}>300/u^{2}$.
Hence \eqref{eq:dUB-from-ass-1a} also follows from the last lower bound for
$d$ above.

Combining the remaining lower bounds for $d$, we have
\begin{equation}
\label{eq:dLB-thm13}
d \geq \max \left(
\frac{0.882 \left| N_{\alpha} \right|^{2/23} b^{48/23}}{u^{44/23}},
\frac{0.66\left| N_{\alpha} \right|^{1/13}b^{28/13}}{u^{24/13}},
\frac{15.3 \left| N_{\alpha} \right|^{4/11}b^{20/11}}{u^{24/11}},
\frac{17 \left| N_{\alpha} \right|^{1/2}b^{2}}{u^{2}}
\right).
\end{equation}

Taking the largest of the coefficients of the terms in this max (which occurs
on the last term), along with the largest exponent on both $b$ and $\left| N_{\alpha} \right|$
in their numerators (occurring in the second and last terms, respectively), as
well as the smallest exponent on $u$ in their denominators (occurring in the
second term) yields
\[
d \geq
\frac{17\left| N_{\alpha} \right|^{1/2} b^{28/13}}{u^{24/13}}.
\]

This is the lower bound for $d$ in Theorem~\ref{thm:1.3-seq-new}(b), so this
completes the proof of that part of the theorem.

\subsection{Proof of Theorem~\ref{thm:1.3-seq-new}(a)}
\label{subsect:proof-thm13a}

As in Subsection~\ref{subsect:thm13-step-v}, we are going to obtain a contradiction
from the assumption that there are three distinct squares
\[
y_{k_{3}}>y_{k_{2}}>y_{k_{1}} \geq \max \left( 4\sqrt{\left| N_{\alpha} \right|/d}, b^{2}\left| N_{\alpha} \right|/d \right)
\]
with $k_{1} \geq 2$ or $k_{1} \leq K-1$. These are the assumptions stated in
\eqref{eq:ass-7aa} and \eqref{eq:ass-5aa} above.

As stated in \eqref{eq:b-LB} at the start of this section, we may assume that
$b \geq 5$. So, with Lemma~\ref{lem:admissible-b}, to complete the proof of
Theorem~\ref{thm:1.3-seq-new}(a), we need to remove the lower bounds on $d$ in
Subsection~\ref{subsect:thm13-step-v} for $b=5,13$ and $17$. We do so by computing.

Rather than the lower bounds for $d$ above in Subsection~\ref{subsect:thm13-step-v},
which depend on $N_{\alpha}$, we will start with lower bounds for $d$ that only
depend on $b$ and $u$.

Instead of \eqref{eq:ass-6a}, \eqref{eq:ass-9a}, \eqref{eq:ass-10a} and \eqref{eq:dLB2-step-iv},
which we used above in Subsection~\ref{subsect:thm13-step-v}, we will use
\eqref{eq:ass-6b}, \eqref{eq:ass-9b}, \eqref{eq:ass-10b} and \eqref{eq:dLB3-step-iv}.
Applying \eqref{eq:b-LB} to \eqref{eq:dUB-from-ass-1a}, we see that
\eqref{eq:dUB-from-ass-1a} holds if $d \geq 12b^{2}/u^{2}$. This inequality
also implies that \eqref{eq:dUB-from-ass-2a} and \eqref{eq:dUB-from-ass-3a}
hold. Hence, we require
\begin{equation}
\label{eq:d-LB2}
d \geq \max \left(
\frac{12b^{2}}{u^{2}},
\frac{0.88 b^{8/3}}{u^{44/21}},
\frac{0.64b^{8/3}}{u^{2}},
\frac{60b^{36/7}}{u^{24/7}}, \frac{85b^{8}}{u^{4}} \right).
\end{equation}

For each value of $b$, we find the smallest value of $u$ such that the right-hand
side of \eqref{eq:d-LB2} is less than $2$ and denote this value of $u$ by $U_{b}$.

\begin{center}
\begin{table}[h]
\begin{tabular}{|c|r|r|r|r|}\hline
$b$ &       $D_{b}$     & $U_{b}$ &   $c_{b}$   & \text{CPU time} \\ \hline
  5 &      $33,203,125$ &   $64$  &    $36,255$ & \text{2s} \\ \hline
 13 &  $69,337,111,285$ &  $432$  & $2,466,430$ & \text{1h 17m 19s} \\ \hline
 17 & $592,939,382,485$ &  $738$  & $7,708,862$ & \text{33h 55m 19s} \\ \hline
\end{tabular}
\caption{Data for $b$}
\label{table:b-data}
\end{table}
\end{center}

We use \eqref{eq:d-LB2} to get initial upper bounds, $D_{b,u}$, on the values
of $d$ that we need to check. For each value of $b$, the maximum value of $D_{b,u}$
comes from $u=1$. We put $D_{b}=D_{b,1}$. The values of $D_{b}$ in Table~\ref{table:b-data}
come from the $85b^{8}/u^{4}$ term in the max in \eqref{eq:d-LB2}.

For each $u$ with $1 \leq u<U_{b}$, consider all integers
$d=\left( t^{2} \pm 4 \right)/ u^{2} \leq D_{b,u}$.

For each such $(b,d,t,u)$, we find all positive integers, $a$, such that $-N_{\alpha}
=db^{4}-a^{2}$ is a positive square -- note that this provides the upper bound $a<\sqrt{db^{4}}$.
We can also use the lower bounds for $d$ in terms of $d$, $N_{\alpha}$ and $u$
in \eqref{eq:ass-6a}, \eqref{eq:ass-9a}, \eqref{eq:ass-10a} and \eqref{eq:dLB2-step-iv}
to get lower bounds for $a$. E.g., from $d>17 \left| N_{\alpha} \right|^{1/2}b^{2}/u^{2}$
in \eqref{eq:dLB2-step-iv}, we obtain $a> \sqrt{db^{4}-\left( du^{2}/ \left( 17b^{2} \right) \right)^{2}}$.
These lower bounds for $a$ significantly reduce the range of values of $a$ that
need to be checked.

We found $c_{b}$ tuples $(a,b,d,t,u)$, where $c_{b}$ is as in Table~\ref{table:b-data}.
For none of these does the lower bound for $d$ in \eqref{eq:dLB-thm13} hold.
Instead we use Proposition~\ref{prop:1.2-seq-new}, which applies to any sequence,
to treat them.

From Proposition~\ref{prop:1.2-seq-new}, we know that for each sequence, there
are at most two distinct squares with $y_{k}$ sufficiently large, where $k \geq 1$
or $k \leq K$. So for each of these tuples, $(a,b,d,t,u)$, for $k \geq 2$ and
$k \leq K-1$, we checked for squares among the $y_{k}$'s not satisfying the lower
bound on $y_{k}$ in Proposition~\ref{prop:1.2-seq-new}. No such sequence were
found.

Initially, the code to search for such tuples $(a,b,d,t,u)$, was written in PARI/GP \cite{Pari}.
However, for $b=13$ and $b=17$, the calculations would have taken too long. So
instead we wrote the code in the Java programming language, which we ran on a
Windows laptop with an Intel i7-13700H CPU and 32~GB of memory.
We record the CPU time for each value of $b$ in Table~\ref{table:b-data}.
The Java code can be found in the file
\verb!BAllNaSqrSearch.java! in the \verb!java! subdirectory of the github url
provided at the end of Section~\ref{sect:intro}.

\end{document}